\documentclass[12pt,reqno]{amsart}

\usepackage[arrow,matrix,curve]{xy}

\usepackage[dvips]{graphicx} 

\usepackage{amssymb, latexsym, amsmath, amscd, array, hyperref
}

\usepackage{setspace,  
}

\theoremstyle{definition}

\numberwithin{equation}{section}

\newcommand\astn {{{}^{\ast}\hspace{-1pt}{\mathbb N}}}

\newcommand\R {{\mathbb R}}

\newcommand\astr {{{}^{\ast}\hspace{-1pt}\mathbb{R}}}

\newcommand\Q {{\mathbb Q}}

\author[J.B.]{Jacques Bair}\address{J. Bair, HEC-ULG, University of
Liege, 4000 Belgium}\email{j.bair@ulg.ac.be}

\author[P.B.]{Piotr B\l{}aszczyk}\address{P. B\l{}aszczyk, Institute
of Mathematics, Pedagogical University of Cracow,
Poland}\email{pb@up.krakow.pl}

\author[K. K.]{Karin U. Katz}\address{K. Katz, Department of
Mathematics, Bar Ilan University, Ramat Gan 52900
Israel}\email{katzmik@math.biu.ac.il}

\author[M.K.]{Mikhail G. Katz}\address{M. Katz, Department of
Mathematics, Bar Ilan University, Ramat Gan 52900 Israel}
\email{katzmik@macs.biu.ac.il}

\author[T.K.]{Taras Kudryk} \address{T. Kudryk, Department of
Mathematics, Lviv National University, Lviv, Ukraine}
\email{kudryk@mail.lviv.ua}

\author[D.S.]{David Sherry}\address{D. Sherry, Department of
Philosophy, Northern Arizona University, Flagstaff, AZ 86011, US}
\email{David.Sherry@nau.edu}

\begin{document}


\pagestyle{empty}  


\title{Analyzing Benardete's comment on decimal notation}

\begin{abstract}
Philosopher Benardete challenged both the conventional wisdom and the
received mathematical treatment of zero, dot, nine recurring.  An
initially puzzling passage in Benardete on the intelligibility of the
continuum reveals challenging insights into number systems, the
foundations of modern analysis, and mathematics education.  A key
concept here is, in Terry Tao's terminology, that of an
\emph{ultralimit}.

Keywords: real analysis; infinitesimals; decimal notation; procedures
\emph{vs} ontology
\end{abstract}

\maketitle


\section{Introduction}

Philosopher Jos\'e Benardete in his book \emph{Infinity: An essay in
metaphysics} argues that some natural pre-mathematical intuitions
cannot be properly expressed if one is limited to an overly
restrictive number system:

\begin{quote} The intelligibility of the continuum has been found--many
times over--to require that the domain of real numbers be enlarged to
include infinitesimals. This enlarged domain may be styled the domain
of continuum numbers. It will now be evident that~$.9999\ldots{}$ does
not equal~$1$ but falls infinitesimally short of it. I think that
$.9999\ldots{}$ should indeed be admitted as a \emph{number} \ldots{}
though not as a \emph{real} number.  \cite[p. 279]{Be64} (emphasis in
the original)
\end{quote}

To a professional mathematician, Benardete's remarks may seem naive.
The equality between~$1$ and~$0.999.\ldots$ is not in the realm of
speculation, but rather an established fact.  This follows from the
very \emph{definition} of~$0.999\ldots$ as the \emph{limit} of the
sequence~$0.9$,~$0.99$,~$0.999, \ldots{}$ where the limit of a
sequence~$(u_n)$ is defined, following any calculus textbook, as the
real number~$L$ such that for every~$\epsilon>0$ there exists an~$N>0$
such that if~$n>N$ then~$|u_n-L|<\epsilon$ (and even in the hyperreal
number system they are still equal, as Bryan Dawson recently noted in
\cite{Da16}).

It patently follows from the definition that~$0.999\ldots{}$ equals
precisely~$1$ because the value~$L=1$ satisfies the definition stated,
and there is nothing more to discuss.  Or is there?  Earlier studies
in this direction include \cite{El10}, \cite{10b}, \cite{10c}.

\section{Vicious circle}

As is patently evident from the definition given above, the procedure
of taking the \emph{limit} is real-valued by definition.  If one is to
treat Benardete's comment charitably, one cannot presuppose the answer
to his query, namely that the number in question is necessarily a real
number.  The equality should, if possible, be conceived independently
from a conception of~$0.999\ldots{}$ in terms of limits that
presupposes that limits are, by definition, real-valued.

While the equality~$0.999\ldots = 1$ (whether it results from a
definition or from a mathematical demonstration) is absolutely
necessary in order to ensure the coherence of the algebra with the
real numbers, in some approaches to the real numbers the
identification of~$0.999\ldots$ with~$1$ is itself part of the
\emph{definition} of the real number line rather than a
\emph{theorem}.  Thus, the real numbers can be defined in terms of
unending strings usually referred to as their \emph{decimal
expansions}.  In this approach to the real numbers, the
strings~$1.000\ldots$ and~$0.999\ldots$ are \emph{postulated} to be
the same number, i.e., they \emph{represent} the same real number;
more generally, the equivalence relation on such formal strings is
defined in such a way that each terminating string is equivalent to
the related one with an infinite tail of~$9s$.

In this approach to the real numbers, it is indeed correct to assert
that the equality (in reality masking an equivalence) is not a theorem
but rather a definition.  It can be debated whether practically
speaking this approach is a good approach to the real numbers;
arguably it is for some purposes, but not for others.  

Historically the first scholar to recognize that it is useful to
represent numbers by \emph{unending} decimals was Simon Stevin already
in the 16th century, before the golden age of the calculus, and even
before the symbolic notation of Vieta; see \cite{KK12b} as well as
\cite[Section~2]{BKS}.  The identification of the two strings is first
found in the second half of the 18th century in the work of Lambert
(1758) and Euler (1765).

\section{Other meanings}

Can one assign any \emph{meaning} to the symbolic string
``$0.999\ldots$'' other that defining it to be~$1$?  That question
cannot be answered without analyzing what \emph{informal} meaning is
assigned to~$0.999\ldots$, \emph{prior to} interpreting it in a formal
mathematical sense.

Beginning calculus students often informally describe this as
``nought, dot, nine recurring'' or alternatively ``zero, point,
followed by infinitely many~$9s$.''  The second description may not
necessarily refer to any sophisticated number system like the real
number system (e.g., equivalence classes of Cauchy sequences of
rational numbers expressed in Zermelo--Fraenkel set theory), since at
this level the students will typically not have been exposed to such
mathematical abstractions, involving as they do equivalence classes of
Cauchy sequences, Dedekind cuts, or other techniques from analysis.

It is also known that at this level, about 80 percent of the students
feel that such an object necessarily falls a little bit short of~$1$.
The question is whether such intuitions are necessarily erroneous, or
whether they could lend themselves to a mathematically rigorous
interpretation in the context of a suitable number system.

A possible mathematical interpretation, in the context of the
sequence~$0.9$,~$0.99$,~$0.999, \ldots{}$, mentioned in the
introduction, is as follows.  Instead of talking the \emph{limit} of
the sequence, one takes, in the terminology of Terry Tao,%
\footnote{See e.g.,
\url{https://terrytao.wordpress.com/tag/ultralimit-analysis}}
an \emph{ultralimit}, resulting in a hyperreal number that falls
infinitesimally short of~$1$.

Rob Ely argues that such intuitions are not necessarily mathematically
erroneous because they can find a rigorous implementation in the
context of a hyperreal number system, where a number with an infinite
tail of 9s can fall infinitesimally short of~$1$ \cite{El10}.  Namely,
if~$H$ is an infinite hypernatural, then the hyperfinite
sum~$\sum_{n=1}^H\frac{9}{10^n} =0.999\ldots9$ contains~$H$
occurrences of the digit~$9$, and falls short of~$1$ by the tiny
amount~$(0.1)^H$.

Infinite hypernaturals like~$H$ belong to a hyperreal line, say
$\astr$, which obeys the same rules, namely the rules of real-closed
ordered fields, as the usual real line~$\R$.  Numbers like~$H$ behave
in~$\astr$ as the usual natural numbers behave in~$\R$, by the
\emph{transfer principle}.%
\footnote{The \emph{transfer principle} is a type of theorem that,
depending on the context, asserts that rules, laws or procedures valid
for a certain number system, still apply (i.e., are ``transfered'') to
an extended number system.  Thus, the familiar
extension~$\Q\hookrightarrow\R$ preserves the properties of an ordered
field.  To give a negative example, the
extension~$\R\hookrightarrow\R\cup\{\pm\infty\}$ of the real numbers
to the so-called \emph{extended reals} does not preserve the
properties of an ordered field.  The hyperreal extension
$\R\hookrightarrow\astr$ preserves \emph{all} first-order properties,
such as the identity~$\sin^2 x + \cos^2 x =1$ (valid for all
hyperreal~$x$, including infinitesimal and infinite values
of~$x\in\astr$).  For a more detailed discussion, see the textbook
\emph{Elementary Calculus} \cite{Ke86}.}
This is in sharp contrast with Cantorian theories where neither
cardinals nor ordinals obey the rules of an ordered field.  Thus,
there is no such thing as~$\omega-1$, while~$\aleph_0+1=\aleph_0$.  On
the other hand, according to Tirosh and Tsamir, students tend to
attribute properties of finite sets to infinite ones:
\begin{quote}
Research in mathematics education indicates that in the transition
from given systems to wider ones learners tend to attribute the
properties that hold for the former to the latter. In particular, it
has been found that, in the context of Cantorian Set Theory, learners
tend to attribute properties of finite sets to infinite ones--using
methods which are acceptable for finite sets, to infinite ones.
\cite{TT}
\end{quote}

The hyperrational number~$\sum_{n=1}^H\frac{9}{10^n} =0.999\ldots9$
can be visualized (virtually) on a hyperreal numerical line in the
\emph{halo} of~$1$ by means of Keisler's microscope, with the
difference being~$10 ^{-H}$.  To return to the real numbers, we have
the equality~$\text{st}(0.999\ldots9) = 1$ where ``st'' is the
operation of taking standard part (shadow).%
\footnote{To elaborate further, note that each hyperreal~$r$ can be
represented by means of its extended decimal expansion~$r=\pm
n.r_1r_2\ldots r_i\ldots r_H\ldots$ where~$n\in\astn\cup\{0\}$
and~$r_i\in\{0,1,2,...,9\}$ for~$i\in\astn$.}

Thus, there are many hyperreal numbers representable as~$0$ followed
by infinitely many~$9$s, e.g.,~$0,999\ldots 9000\ldots$ (where the
zeros start from some infinite rank),~$0.999\ldots 9123123\ldots$,
$0,999\ldots 999\ldots$, etc.  Among such numbers only one is the
natural extension of the usual sequence~$0.(9)$ (unending~$9$s) and
this indeed equals~$1$.  The matter was explored in more detail by
\cite{Li}.

The hyperrational number~$\sum_{n=1}^H\frac{9}{10^n} =0.999\ldots9$ is
a \emph{terminating} infinite string of~$9$s different from the one
usually envisioned in real analysis, but it respects student
intuitions.  Intuitions like these underlie Leibniz's and Peirce's
conceptions of the continuum.  They can be helpful in learning the
calculus, as argued in \cite{El10}; see also \cite{Vi16}, \cite{KP}.

The pedagogical issue is a separate one but what could be emphasized
here is that the existence of such an interpretation suggests that we
indeed do \emph{assume} that such a string represents a \emph{real}
number when we prove that it necessarily equals~$1$.  The idea that
such an assumption can be challenged is in line with Benardete's
comment.

\section{Conclusion}

It is known that the real number system provides an adequate
foundation for mathematics.  Yet the real number system does not
explain the well-known phenomenon of student unease about the identity
$1=0.999\ldots$ (nor the philosopher's unease).

It is similarly known that a set-theoretic justification of a
hyperreal number system involves more work than is justified by the
need of explaining~$0.999\ldots{}$ to the students.  Yet the crucial
distinction between set-theoretic justification (ontology), on the one
hand, and the \emph{procedures} involving infinitesimals, on the
other, has been emphasized in \cite{Ba16} and \cite{Ba16b}.

Working with infinitesimals and exploiting the procedures based on
them to great scientific effect was a salient feature of the output of
the pioneers of analysis like Fermat (see \cite{13e}), Leibniz (see
\cite{14a}), Euler (see \cite{13a}), and Cauchy (see \cite{18a}).

What we show in our paper is that such infinitesimal procedures can be
used to address both the student unease and the philosopher's unease.
Leibniz did not develop a set-theoretic ontology for infinitesimals
for the simple reason that set theory did not exist yet, but he may
have been more open to an infinitesimal gap between ``zero, dot,
followed by an infinite number of~$9$s'' and~$1$ than many
traditionally trained mathematicians today.

\end{document}